\documentclass{amsart}%
\usepackage{amsthm,amstext,amscd,latexsym,amsmath,amsfonts,amssymb}
\usepackage{amsmath}
\usepackage{amsfonts}
\usepackage{amssymb}%

\theoremstyle{plain}
\newtheorem{theorem}{Theorem}
\newtheorem{lemma}[theorem]{Lemma}
\newtheorem{proposition}[theorem]{Proposition}

\theoremstyle{definition}
\newtheorem{definition}[theorem]{Definition}
\newtheorem{corollary}[theorem]{Corollary}
\newtheorem{example}{\sc Example}
\theoremstyle{remark}
\newtheorem{remark}{\sc Remark}

\begin{document}
\title{Bi-element representations of ternary groups}
\author{Andrzej Borowiec}
\address{Institute of Theoretical Physics, University of Wroc{\l }aw,
Pl. Maxa Borna 9,
50-204 Wroc{\l }aw, Poland}
\email{borow@ift.univ.wroc.pl}
\author{Wies{\l }aw A. Dudek}
\address{Institute of Mathematics, Technical University of Wroc{\l }aw, Wybrzeze
Wyspianskiego 27, 50-370 Wroc{\l }aw, Poland}
\email{dudek@im.pwr.wroc.pl}
\author{Steven Duplij}
\address{Department of Physics and Technology, Kharkov National University, Kharkov
61001, Ukraine}
\email{Steven.A.Duplij@univer.kharkov.ua}
\urladdr{http://www.math.uni-mannheim.de/\symbol{126}duplij}

\maketitle
\begin{quotation}
\begin{center} {\it
Dedicated to the memory of our coauthor and friend \\W\l adek Marcinek 
who unexpectedly passed away on June 9, 2003 }
\end{center}\end{quotation}
\mbox{}
\medskip

\begin{abstract}
General properties of ternary semigroups and groups are considered. The
bi-element representation theory in which every representation matrix
corresponds to a pair of elements is built, connection with the standard
theory is considered and several concrete examples are constructed. For
clarity the shortened versions of classical Gluskin-Hossz\'{u} and Post
theorems are given for them.
\end{abstract}

\section{Introduction}

Ternary and $n$-ary generalizations of algebraic structures is the most
natural way for further development and deeper understanding of their
fundamental properties. Firstly ternary algebraic operations were introduced
already in the XIX-th century by A. Cayley. As the development of Cayley's
ideas it were considered $n$-ary generalization of matrices and their
determinants \cite{sokolov,kap/gel/zel} and general theory of $n$-ary algebras
\cite{law,car5} and ternary rings \cite{lis}. For some physical applications in
Nambu mechanics, supersymmetry, Yang-Baxter equation, etc. see  e.g.
\cite{ker1,vai/ker,ozi/paa/roz}. The statement here is something
different and based on
our previous preliminary report \cite{bor/dud/dup1}, where also
ternary algebras and ternary Hopf algebras were considered.

The notion of an $n$-ary group
was introduced in 1928 by W. D\"{o}rnte \cite{dor3} (inspired by E.
N\"{o}ther) which is a natural generalization of the notion of a group and a
ternary group considered by Certaine \cite{cer} and Kasner \cite{kas}.

In this paper we reformulate necessary material on ternary semigroups and
groups \cite{belousov,rusakov1} in abstract language.

In the (binary) group theory representations were introduced as a matrix
realization of group elements and abstract group action by usual matrix
multiplication, when one element was described by one matrix. In this paper we
propose a new bi-element approach to the representation theory of ternary
group, when one matrix parametrizes two elements of a ternary group.

An alternative approach to ternary group representations was made in
\cite{wan/wan}, where it was proposed, instead of our operator-valued
functions of two variables, functions of one variable taking value in a pair
of operators (matrices), viz. $\Pi^{wan}:G\rightarrow\left(  \Pi_{1}%
^{wan}\left(  x\right)  ,\Pi_{2}^{wan}\left(  x\right)  \right)  \in GL\left(
V\right)  \times GL\left(  V\right)  $ with another analog of homomorphism.
Unfortunately, in \cite{wan/wan} there were no given concrete examples nor
connections with derived case.

Here, using our method, we present several concrete examples and consider
connection with binary case proving the classical Gluskin-Hossz\'{u}
and Post theorems.

\section{Ternary semigroups}

A non-empty set $G$ with one \emph{ternary} operation $[\;]:G\times G\times
G\rightarrow G$ is called a \emph{ternary groupoid} and is denoted by
$(G,[\;])$ or $\left(  G,m^{\left(  3\right)  }\right)  $. We will present
some results using second notation, because it allows to reverse arrows in the
most clear way. In proofs we will mostly use the first notation due to
convenience and for short.

If on $G$ there is a binary operation $\odot$ (i.e. $m^{\left(  2\right)  }$)
such that $[xyz]=(x\odot y)\odot z$, i.e.%
\begin{equation}
m^{\left(  3\right)  }=m^{\left(  2\right)  }\circ\left(  m^{\left(  2\right)
}\times\operatorname*{id}\right)  \label{mder}%
\end{equation}
for all $x,y,z\in G$, then we say that $[\;]$ (respectively $m^{\left(
3\right)  }$) is \emph{derived} from $\odot$ (respectively from $m^{\left(
2\right)  })$ and denote this fact by $(G,[\;])=\mathrm{der}(G,\odot)$
(respectively by $m^{\left(  3\right)  }=m_{\mathrm{der}}^{\left(  3\right)
}$). If
\[
\lbrack xyz]=((x\odot y)\odot z)\odot b
\]
holds for all $x,y,z\in G$ and some fixed $b\in G$, then a groupoid $(G,[\;]$
is \emph{$b$-derived} from $(G,\odot)$. In this case we write
$(G,[\;])=\mathrm{der}_{b}(G,\odot)$ (cf. \cite{dud/mic1,dud/mic2}).
%{DM1)
%W.A.Dudek, J.Michalski: On a generalization of Hoss\' theorem,
%Demonstratio Math. {\bf 15} (1982), 437-441.
%{DM2)
%W.A.Dudek, J.Michalski: On retract of polyadic groups,
%Demonstratio Math. {\bf 17} (1984), 281-301.

We say that $(G,[\;]$ is a \emph{ternary semigroup} if the operation $[\;]$ is
\emph{associative}, i.e. if
\begin{equation}
\left[  \left[  xyz\right]  uv\right]  =\left[  x\left[  yzu\right]  v\right]
=\left[  xy\left[  zuv\right]  \right]  \label{ass1}%
\end{equation}
holds for all $x,y,z,u,v\in G$, i.e.%
\begin{equation}
m^{\left(  3\right)  }\circ\left(  m^{\left(  3\right)  }\times
\operatorname*{id}\times\operatorname*{id}\right)  =m^{\left(  3\right)
}\circ\left(  \operatorname*{id}\times m^{\left(  3\right)  }\times
\operatorname*{id}\right)  =m^{\left(  3\right)  }\circ\left(
\operatorname*{id}\times\operatorname*{id}\times m^{\left(  3\right)
}\right)  . \label{assm}%
\end{equation}

Obviously, a ternary operation $m_{\mathrm{der}}^{\left(  3\right)  }$ derived
from a binary associative operation $m^{\left(  2\right)  }$ is also
associative in the above sense, but a ternary operation $[\;]$ which is
$b$-derived from an associative operation $\odot$ is associative in the above
sense, if and only if $b$ lies in the center of $(G,\odot)$.

Fixing one element in a ternary operation we obtain a binary operation. A
binary groupoid $(G,\odot)$, where $x\odot y=[xay]$ for some fixed $a\in G$,
respectively $\left(  G,m_{a}^{\left(  2\right)  }\right)  $, where
\begin{equation}
m_{a}^{\left(  2\right)  }=m^{\left(  3\right)  }\circ\left(
\operatorname*{id}\times a\times\operatorname*{id}\right)  ,\label{ma}%
\end{equation}
is called a \emph{retract} of $(G,[\;])$ and is denoted by $\mathrm{ret}%
_{a}(G,[\;])$. In some special cases described in \cite{dud/mic1,dud/mic2} we
have $(G,\odot)=\mathrm{ret}_{a}(\mathrm{der}_{b}(G,\odot))$ and
$(G,\odot)=\mathrm{ret}_{c}(\mathrm{der}_{d}(G,\odot))$, but in general
$(G,\odot)$ and $\mathrm{ret}_{a}(\mathrm{der}_{b}(G,\odot))$ are only
isomorphic \cite{dud/mic2}.

\begin{lemma}
\label{lem1}If in the ternary semigroup $(G,[\;])$ there exists an element $e$
such that for all $y\in G$ we have $\left[  eye\right]  =y$, then this
semigroup is derived from the binary semigroup $\mathrm{ret}_{e}(G,[\;])$,
i.e. $(G,[\;])=\mathrm{der}(\mathrm{ret}_{e}(G,[\;]))$, and this semigroup is
derived from the binary semigroup $(G,m_{e}^{\left(  2\right)  })$, where%
\begin{equation}
m_{e}^{\left(  2\right)  }=m^{\left(  3\right)  }\circ\left(
\operatorname*{id}\times e\times\operatorname*{id}\right)  .\label{me}%
\end{equation}

\end{lemma}

\begin{proof}
Indeed, if we put $x\circledast y=[xey]$, then $(x\circledast y)\circledast
z=[[xey]ez]=[x[eye]z]=[xyz]$ and $x\circledast(y\circledast
z)=[xe[yez]]=[x[eye]z]=[xyz]$, which completes the proof.
\end{proof}

The same ternary semigroup $\left(  G,m^{\left(  3\right)  }\right)  $ can be
derived from two different (but isomorphic) semigroups $\left(  G,\circledast
\right)  $ and $\left(  G,\diamond\right)  $ ($\left(  G,m_{e}^{\left(
2\right)  }\right)  $ and $\left(  G,m_{a}^{\left(  2\right)  }\right)  $).
Indeed, if in $G$ there exists $a\neq e$ such that $[aya]=y$ for all $y\in G$,
then by the same argumentation we obtain $[xyz]=x\diamond y\diamond z$ for
$x\diamond y=[xay]$. In this case for $\varphi(x)=x\diamond e=[xae]$ we have
\[
x\circledast y=[xey]=[x[aea]y]=[[xae]ay]=(x\diamond e)\diamond y=\varphi
(x)\diamond y
\]
and
\[
\varphi(x\circledast y)=[[xey]ae]=[[x[aea]y]ae]=[[xae]a[yae]]=\varphi
(x)\diamond\varphi(y).
\]
Thus $\varphi$ is a binary homomorphism such that $\varphi(e)=a$. Moreover,
for $\psi(x)=[eax]$ we have
\begin{align*}
\psi(\varphi(x)) &  =[ea[xae]]=[e[axa]e]=x,\\
\varphi(\psi(x)) &  =[[eax]ae]=[e[axa]e]=x
\end{align*}
and
\[
\psi(x\diamond y)=[ea[xay]]=[ea[x[eae]y]]=[[eax]e[aey]]=\psi(x)\circledast
\psi(y).
\]
Hence semigroups $(G,\circledast)$ and $(G,\diamond)$ are isomorphic.

\begin{definition}
An element $e\in G$ is called a \emph{middle identity} or a \emph{middle
neutral element} of $(G,[\;])$, if for all $x\in G$ we have $[exe]=x$, i.e.%
\begin{equation}
m^{\left(  3\right)  }\circ\left(  e\times\operatorname*{id}\times e\right)
=\operatorname*{id}.\label{mmx}%
\end{equation}
An element $e\in G$ satisfying the identity $[eex]=x$, i.e.%
\begin{equation}
m^{\left(  3\right)  }\circ\left(  e\times e\times\operatorname*{id}\right)
=\operatorname*{id}\label{ml}%
\end{equation}
is called a \emph{left identity} or a \emph{left neutral element} of
$(G,[\;])$. Similarly we define a \emph{right identity}. An element which is a
left, middle and right identity is called a \emph{ternary identity} (briefly: identity).
\end{definition}

There are ternary semigroups without left (middle, right) neutral elements,
but there are also ternary semigroups in which all elements are identities.

\begin{example}
In ternary semigroups derived from the symmetric group $S_{3}$ all elements of
order 2 are left and right (but no middle) identities.
\end{example}

\begin{example}
In ternary semigroup derived from Boolean group all elements are ternary
identities, but ternary semigroup $1$-derived from the additive group
$\mathbb{Z}_{4}$ has no left (right, middle) identities.
\end{example}

\begin{lemma}
For any ternary semigroup $(G,[\;])$ with a left (right) identity there exists
a binary semigroup $(G,\odot)$ and its endomorphism $\mu$ such that
\[
\lbrack xyz]=x\odot\mu(y)\odot z
\]
for all $x,y,z\in G$.
\end{lemma}

\begin{proof}
Let $e$ be a left identity of $(G,[\;])$. It is not difficult to see that the
operation $x\odot y=[xey]$ is associative. Moreover, for $\mu(x)=[exe]$, we
have
\[
\mu(x)\odot\mu(y)=[[exe]e[eye]]=[[exe][eey]e]=[e[xey]e]=\mu(x\odot y)
\]
and
\[
\lbrack xyz]=[x[eey][eez]]=[[xe[eye]]ez]=x\odot\mu(y)\odot z.
\]

In the case of right identity the proof is analogous.
\end{proof}

\begin{definition}
We say that a ternary groupoid $(G, [\; ])$ is:

a \emph{left cancellative} if $[abx] = [aby]\Longrightarrow x=y$,

a \emph{middle cancellative} if $[axb] = [ayb]\Longrightarrow x=y$,

a \emph{right cancellative} if $[xab] = [yab]\Longrightarrow x=y$\newline
holds for all $a,b\in G$.

A ternary groupoid which is left, middle and right cancellative is called
\emph{cancellative}.
\end{definition}

\begin{theorem}
\label{theor-canc}A ternary groupoid is cancellative if and only if it is a
middle cancellative, or equivalently, if and only if it is a left and right cancellative.
\end{theorem}

\begin{proof}
Assume that a ternary semigroup $(G,[\;])$ is a middle cancellative and
$[xab]=[yab]$. Then $[ab[xab]]=[ab[yab]]$ and in the consequence
$[a[bxa]b]=[a[bya]b]$ which implies $x=y$.

Conversely if $(G,[\;])$ is a left and right cancellative and $[axb]=[ayb]$
then $[a[axb]b]=[a[ayb]b]$ and $[[aax]bb]=[[aay]bb]$ which gives $x=y$.
\end{proof}

The above theorem is a consequence of the general result proved in
\cite{dud1}.
%
%W.A. Dudek: Autodistributive $n$-groups, Annales Sci. math. Polonae,
%Commentationes Math. 23 (1983), 1-11.

\begin{definition}
A ternary groupoid $(G,[\ ])$ is called $\sigma$-commutative, if
\begin{equation}
\left[  x_{1}x_{2}x_{3}\right]  =\left[  x_{\sigma\left(  1\right)  }%
x_{\sigma\left(  2\right)  }x_{\sigma\left(  3\right)  }\right]  \label{xxx}%
\end{equation}
holds for all $x_{1},x_{2},x_{3}\in G$, i.e. if $m^{(3)}=m^{(3)}\circ\sigma$.
If (\ref{xxx}) holds for all $\sigma\in S_{3}$, then $(G,[\ ])$ is a
\emph{commutative} groupoid. If (\ref{xxx}) holds only for $\sigma=(13)$, i.e.
if $[x_{1}x_{2}x_{3}]=[x_{3}x_{2}x_{1}]$, then $(G,[\ ])$ is called
\emph{semicommutative}.
\end{definition}

The group $S_{3}$ is generated by two transpositions; $(12)$ and $(23)$. This
means that $(G,[\;])$ is commutative if and only if $[xyz]=[yxz]=[xzy]$ holds
for all $x,y,z\in G$.

As a simple consequence of Theorem \ref{theor-canc} from \cite{dud2}
%
%Rem
%w.A. Dudek: Remarks on $n$-groups, Demonstratio Math. 13 (1980), 165-181.
we obtain

\begin{corollary}
If in a ternary semigroup $(G,[\; ])$ satisfying the identity $[xyz]=[yxz]$
there are $a,b$ such that $[axb]=x$ for all $x\in G$, then $(G,[\; ])$ is commutative.
\end{corollary}

\begin{proof}
According to the above remark it is sufficient to prove that $[xyz]=[xzy]$. We
have
\[
\lbrack xyz]=[a[xyz]b]=[ax[yzb]]=[ax[zyb]]=[a[xzy]b]=[xzy].
\]

\end{proof}

Mediality in the binary case is%
\begin{equation}
\left(  x\odot y\right)  \odot\left(  z\odot u\right)  =\left(  x\odot
z\right)  \odot\left(  y\odot u\right)  . \label{med}%
\end{equation}
This can be presented as a matrix $A^{\left(  2\right)  }=\left(
\begin{array}
[c]{cc}%
x & y\\
z & u
\end{array}
\right)  $, read from left by rows and from top by columns as $%
\begin{array}
[c]{ccc}
& \Downarrow & \Downarrow\\
\Rightarrow & x & y\\
\Rightarrow & z & u
\end{array}
$ (see \cite{belousov}).

\begin{remark}
In the binary case a middle cancellative semigroup is commutative, and so for
groups mediality coincides with commutativity.
\end{remark}

In the ternary case instead of $A^{\left(  2\right)  }$ we have $3\times3$
matrix $A^{\left(  3\right)  }$ which should be read similarly.

\begin{definition}
A ternary groupoid $(G,[\;])$ is \emph{medial} if it satisfies the identity
\[
\lbrack\lbrack x_{11}x_{12}x_{13}][x_{21}x_{22}x_{23}][x_{31}x_{32}%
x_{33}]]=[[x_{11}x_{21}x_{31}][x_{12}x_{22}x_{32}][x_{13}x_{23}x_{33}]],
\]
i.e.%
\begin{equation}
m^{\left(  3\right)  }\circ\left(  m^{\left(  3\right)  }\times m^{\left(
3\right)  }\times m^{\left(  3\right)  }\right)  =m^{\left(  3\right)  }%
\circ\left(  m^{\left(  3\right)  }\times m^{\left(  3\right)  }\times
m^{\left(  3\right)  }\right)  \circ\sigma_{medial}, \label{smm}%
\end{equation}
where $\sigma_{medial}=\binom{123456789}{147258369}\in S_{9}.$
\end{definition}

It is not difficult to see that a semicommutative ternary semigroup is medial.

An element $x$ such that $[xxx]=x$ is called an \emph{idempotent}. A groupoid
in which all elements are idempotents is called an \emph{idempotent groupoid}.
A left (right, middle) identity is an idempotent.

\section{Ternary groups}

\begin{definition}
A ternary semigroup $(G,[\;])$ is a \emph{ternary group} if for all $a,b,c\in
G$ there are $x,y,z\in G$ such that
\begin{equation}
\lbrack xab]=[ayb]=[abz]=c. \label{ac}%
\end{equation}

\end{definition}

One can prove \cite{pos}
%
%post
%E.L.Post: Polyadic groups, Trans. amer. Math. soc. 48 (1940), 208-350.
that elements $x,y,z$ are uniquely determined. Moreover, according to the
suggestion of \cite{pos} one can prove (cf. \cite{dud/gla/gle})
%
%DGG
%W.A.Dudek, K.G{\l}azek, B.Gleichgewicht: A note on the axioms of $n$-groups,
%Coll. Math. soc. J. Bolyai. 29. Universal Algebra, Esztergom (Hungary) 1977,
%pp. 195-202.
that in the above definition, under the assumption of the associativity, it
suffices only to postulate the existence of a solution of $[ayb]=c$, or
equivalently, of $[xab]=[abz]=c$.

In a ternary group the equation $[xxz]=x$ has a unique solution which is
denoted by $z=\overline{x}$ and called \emph{skew element} (cf. \cite{dor3}),
or in the other notation%
\[
m^{\left(  3\right)  }\circ\left(  \operatorname*{id}\times\operatorname*{id}%
\times\overline{\cdot}\right)  \circ D^{\left(  3\right)  }=\operatorname*{id}%
,
\]
where $D^{\left(  3\right)  }\left(  x\right)  =\left(  x,x,x\right)  $ is a
ternary diagonal map and $\overline{\cdot}:x\rightarrow\bar{x}$. As a
consequence of results obtained in \cite{dor3} we have

\begin{theorem}
In any ternary group $(G,[\;])$ for all $x,y,z\in G$ the following relations
take place
\begin{align*}
\lbrack x\,x\,\overline{x}]  &  =[x\,\overline{x}\,x]=[\,\overline
{x}\,x\ x]=x,\\
\lbrack y\,x\,\overline{x}]  &  =[y\,\overline{x}\,x]=[x\,\overline
{x}\,y]=[\,\overline{x}\,x\,y]=y,\\
\overline{\lbrack x\ y\ z]}  &  =[\,\overline{z}\ \overline{y}\ \overline
{x}],\\
\overline{\overline{x}}  &  =x.
\end{align*}

\end{theorem}

Since in an idempotent ternary group $\overline{x}=x$ for all $x$, an
idempotent ternary group is semicommutative. From the results obtained in
\cite{dud/gla/gle} (see also \cite{dud2}) for $n=3$ we have

\begin{theorem}
A ternary semigroup $(G,[\;])$ with a unary operation ${}^{-}:x\rightarrow
\overline{x}$ is a ternary group if and only if it satisfies identities
\[
\lbrack yx\,\overline{x}\,]=[x\,\overline{x}\,y]=y,
\]
or in other notation%
\begin{align*}
m^{\left(  3\right)  }\circ\left(  \operatorname*{id}\times\operatorname*{id}%
\times\overline{\cdot}\right)  \circ\left(  \operatorname*{id}\times
D^{\left(  2\right)  }\right)   &  =\Pr\nolimits_{1},\\
m^{\left(  3\right)  }\circ\left(  \operatorname*{id}\times\overline{\cdot
}\times\operatorname*{id}\right)  \circ\left(  D^{\left(  2\right)  }%
\times\operatorname*{id}\right)   &  =\Pr\nolimits_{2},
\end{align*}
where $D^{\left(  2\right)  }\left(  x\right)  =\left(  x,x\right)  $ and
$\Pr\nolimits_{1}\left(  x,y\right)  =x,$ $\Pr\nolimits_{2}\left(  x,y\right)
=y$.
\end{theorem}

\begin{corollary}
A ternary semigroup $(G,[\;])$ is an idempotent ternary group if and only if
it satisfies identities
\[
\lbrack yxx]=[xxy]=y.
\]

\end{corollary}

By Lemma \ref{lem1} a ternary group with an identity is derived from a binary group.

\begin{remark}
The set $S_{3}\backslash A_{3}$ of all odd permutations with ternary operation
$[\;]$ defined as composition of three permutations is an example of a
noncommutative ternary group which is not derived from any group (all groups
with three elements are commutative and isomorphic to $\mathbb{Z}_{3}$).
\end{remark}

From results proved in \cite{dud2} follows

\begin{theorem}
A ternary group $(G,[\;])$ satisfying the identity
\[
\lbrack xy\overline{x}]=y
\]
or
\[
\lbrack\overline{x}yx]=y
\]
is commutative.
\end{theorem}

The most important theorem is

\begin{theorem}
[Gluskin-Hossz\'{u}]\label{thm-gh}For a ternary group $\left(  G,\left[
\;\right]  \right)  $ and fixed element $a\in G$ there exist a binary group
$\left(  G,\circledast\right)  =\mathrm{ret}_{a}\left(  G,\left[  \;\right]
\right)  $ and its automorphism $\varphi$ such that $\varphi(a)=a$ and
\begin{equation}
\left[  xyz\right]  =x\circledast\varphi\left(  y\right)  \circledast
\varphi^{2}\left(  z\right)  \circledast b,\label{gh}%
\end{equation}
where $b=[\overline{a}\,\overline{a}\,\overline{a}\,]$.
\end{theorem}

\begin{proof}
Let $a\in G$ be fixed. Then the binary operation $x\circledast y=\left[
x\,a\,y\right]  $ is associative, because
\[
(x\circledast y)\circledast z=[[xay]az]=[xa[yaz]]=x\circledast(y\circledast
z).
\]
In $(G,\circledast)$ an element $\overline{a}$ is the identity, $[\overline
{a}\,\overline{x}\,\overline{a}]$ inverse of $x$. $\varphi(x)=\left[
\overline{a}\ x\ a\right]  $ is an automorphism of $(G,\circledast)$. The easy
calculation proves that the above formula holds for $b=[\overline
{a}\,\overline{a}\,\overline{a}\,]$. (see \cite{sok} and \cite{dud/mic1}).
\end{proof}

%Sok
%E.I.Sokolov: On the theorem of Gluskin-Hossz\'u on D\"ornte groups,
%(Russian), Mat. Issled. 39 (1976), 187-189.
One can prove that the group $(G,\circledast)$ is unique up to isomorphism
\cite{dud/mic1}. From the proof of Theorem 3 in \cite{gla/gle}
%
%GG
%K.G{\l}azek, B. Gleichgewicht: Abelian $n$-groups,  Coll. math. soc. j. Bolyai,
%29. Universal algebra, Esztergom (Hungary) 1977, 321-329.
it follows that any medial ternary group satisfies the identity
\[
\overline{\,[xyz]\,}=[\overline{x}\,\overline{y}\,\overline{z}\,],
\]
which together with our previous results shows that in such groups we have
\[
\lbrack\overline{x}\,\overline{y}\,\overline{z}]=[\overline{z}\,\overline
{y}\,\overline{x}].
\]
But $\overline{\overline{x}}=x$. Hence, any medial ternary group is
semicommutative, thus any retract of such group is a commutative group.
Moreover, for $\varphi$ from the proof of Theorem \ref{thm-gh} and
$\varphi(b)=b$ for $b=[\overline{a}\,\overline{a}\,\overline{a}]$ we have
\[
\varphi(\varphi(x))=[\overline{a}\,[\overline{a}xa]a]=[\overline
{a}\,a\,[x\overline{a}\,a]]=x
\]

\begin{corollary}
Any medial ternary group $(G,[\;])$ has the form
\[
\lbrack xyz]=x\odot\varphi(y)\odot z\odot b,
\]
where $(G,\odot)$ is a commutative group, $\varphi$ its automorphism such that
$\varphi^{2}=\operatorname*{id}$ and $b\in G$ is fixed.
\end{corollary}

\begin{corollary}
A ternary group is medial, if and only if it is semicommutative.
\end{corollary}

\begin{corollary}
\label{medial}A ternary group is semicommutative (medial), if and only if
there exists $a\in G$ such that $[xay]=[yax]$ holds for all $x,y\in G$.
\end{corollary}

\begin{corollary}
A commutative ternary group is $b$-derived from some commutative group.
\end{corollary}

Indeed, $\varphi(x)=[\overline{a}\,xa]=[xa\,\overline{a}]=x$.

\begin{theorem}
[Post]\label{theor-post}For any ternary group $\left(  G,\left[  \ \right]
\right)  $ there exists a binary group $\left(  G^{\ast},\circledast\right)  $
and $H\lhd G^{\ast}$, such that $G^{\ast}\diagup H\simeq\mathbb{Z}_{2}$ and
\[
\left[  xyz\right]  =x\circledast y\circledast z
\]
for all $x,y,z\in G$.
\end{theorem}

\begin{proof}
Let $c$ be a fixed element in $G$ and let $G^{\ast}=G\times\mathbb{Z}_{2}$. In
$G^{\ast}$ we define binary operation $\circledast$ putting
\[
(x,0)\circledast(y,0) = ([xy\overline{c}], 1 )
\]
\[
(x,0)\circledast(y,1) = ([xyc], 0 )
\]
\[
(x,1)\circledast(y,0) = ([xcy], 0 )
\]
\[
(x,1)\circledast(y,1) = ([xcy], 1 ).
\]

It is not difficult to see that this operation is associative and
$(\overline{c},1)$ is its neutral element. The inverse element (in $G^{\ast}$)
has the form:
\[
(x,0)^{-1} = (\overline{x},0)
\]
\[
(x,1)^{-1} = ([\overline{c}\,\overline{x}\,\overline{c}], 1)
\]

Thus $G^{\ast}$ is a group such that $H=\{(x,1):x\in G\}\vartriangleleft
G^{\ast}$. Obviously the set $G$ can be identified with $G\times\{0\}$ and
\begin{align*}
x\circledast y\circledast z  &  =((x,0)\circledast(y,0))\circledast
(z,0)=([xy\overline{c}],1)\circledast(z,0)\\
&  =([[xy\overline{c}]cz],0)=([xy[\overline{c}cz]],0)=([xyz],0)=[xyz],
\end{align*}
which completes the proof.
\end{proof}

The original proof of this theorem uses some equivalences of sequences of
elements from $G$ (see \cite{pos}). Our proof is based on some general method
presented in \cite{mic1}. Note that group $G^{\ast}$ satisfying all conditions
formulated in our theorem is called \emph{covering} for ternary group
$(G,\left[  \ \right]  )$. Our construction gives the free covering group in
this sense of universal algebras. From results obtained in \cite{dud/mic2} it follows

\begin{proposition}
All retracts of a ternary group $(G,[\ ])$ are isomorphic to the normal
subgroup $H$ of $G^{\ast}$ from the previous theorem, i.e.
\[
\mathrm{ret}_{a}\left(  G,\left[  \;\right]  \right)  \simeq H\vartriangleleft
G^{\ast}.
\]

\end{proposition}

\section{Binary representations of ternary groups}

For a given ternary group $\left(  G,\left[  \;\right]  \right)  $ denote by
$\left(  G\times G,\ast\right)  $ a semigroup with the following binary
multiplication%
\begin{equation}
\left(  x,y\right)  \ast\left(  u,v\right)  =\left(  \left[  xyu\right]
,v\right)  .\label{gg}%
\end{equation}
Obviously, for all $x,u,v\in G$ we have $\left(  x,\overline{x}\right)
\ast\left(  u,v\right)  =\left(  \overline{x},x\right)  \ast\left(
u,v\right)  =\left(  u,v\right)  $, which means that $\left(  x,\overline
{x}\right)  $ and $\left(  \overline{x},x\right)  $ are left (but not right)
unities in $\left(  G\times G,\ast\right)  $. Generally $\left(
x,\overline{x}\right)  \neq\left(  \overline{x},x\right)  $. But for all
$x,y\in G$ we have also $\left(  x,y\right)  \ast\left(  \overline
{y},y\right)  =\left(  \overline{y},y\right)  \ast\left(  x,y\right)  =\left(
x,y\right)  $, i.e. each element $\left(  x,y\right)  $ has a "private" unit.
Moreover, any element $\left(  u,\bar{u}\right)  $, $u\in G$ is a left unit.

The semigroup $\left(  G\times G,\ast\right)  $ is left (but not right)
cancellative, i.e. $\left(  a,b\right)  \ast\left(  x,y\right)  =\left(
a,b\right)  \ast\left(  c,d\right)  $ implies $\left(  x,y\right)  =\left(
c,d\right)  $. Moreover, $\left(  G\times G,\ast\right)  $ is also a right
quasigroup, i.e. for every $\left(  a,b\right)  ,\left(  c,d\right)  \in
G\times G$ there exists only one $\left(  x,y\right)  \in G\times G$ such that
$\left(  a,b\right)  \ast\left(  x,y\right)  =\left(  c,d\right)  $. Similarly
it is not difficult to see that for each $a,b,c,d\in G$ there are uniquely
determined $x,y\in G$ such that $\left(  x,a\right)  \ast\left(  b,c\right)
=\left(  a,y\right)  \ast\left(  b,c\right)  =\left(  d,c\right)  $.

Let $V$ be a vector space over $\mathbb{K}$ and $\operatorname*{End}V$ be a
set of linear endomorphisms of $V$.

\begin{definition}
A \textit{left representation} of a ternary group $(G,[\ ])$ in $V$ is a map
$\Pi^{L}:G\times G\rightarrow\operatorname*{End}V$ such that%
\begin{align}
\Pi^{L}\left(  x_{1},x_{2}\right)  \circ\Pi^{L}\left(  x_{3},x_{4}\right)   &
=\Pi^{L}\left(  \left[  x_{1}x_{2}x_{3}\right]  ,x_{4}\right)  ,\;\;\forall
x_{1},x_{2},x_{3},x_{4}\in G\label{p1}\\
\Pi^{L}\left(  x,\overline{x}\right)   &  =\operatorname*{id}\nolimits_{V}%
,\;\;\forall x\in G.\label{p2}%
\end{align}

\end{definition}

Replacing in (\ref{p2}) $x$ by $\overline{x}$ we obtain $\Pi^{L}\left(
\overline{x},x\right)  =id_{V}$, which means that in fact (\ref{p2}) has the
form $\Pi^{L}\left(  \overline{x},x\right)  =\Pi^{L}\left(  x,\overline
{x}\right)  =\operatorname*{id}\nolimits_{V},\;\;\forall x\in G.$

\begin{lemma}
\label{xxl}For all $x_{1},x_{2},x_{3},x_{4}\in G$ we have
\[
\Pi^{L}\left(  \left[  x_{1}x_{2}x_{3}\right]  ,x_{4}\right)  =\Pi^{L}\left(
x_{1},\left[  x_{2}x_{3}x_{4}\right]  \right)  .
\]

\end{lemma}

\begin{proof}
Indeed, we have%
\begin{align*}
\Pi^{L}\left(  \left[  x_{1}x_{2}x_{3}\right]  \,,x_{4}\right)   &  =\Pi
^{L}\left(  \left[  x_{1}x_{2}x_{3}\right]  \,,x_{4}\right)  \circ\Pi
^{L}\left(  x,\overline{x}\right) \\
&  =\Pi^{L}\left(  \left[  \left[  x_{1}x_{2}x_{3}\right]  x_{4}\,x\right]
,\,\overline{x}\right)  =\Pi^{L}\left(  \left[  x_{1}\left[  x_{2}x_{3}%
x_{4}\right]  x\right]  ,\,\overline{x}\right) \\
&  =\Pi^{L}\left(  x_{1},\left[  x_{2}x_{3}x_{4}\right]  \right)  \circ\Pi
^{L}\left(  x,\overline{x}\right)  =\Pi^{L}\left(  x_{1},\left[  x_{2}%
x_{3}x_{4}\right]  \right)  .
\end{align*}

\end{proof}

Note also that for all $x,y,z\in G$ we have%

\begin{equation}
\Pi^{L}\left(  x,y\right)  =\Pi^{L}\left(  [xz\overline{z}],y\right)  =\Pi
^{L}\left(  x,z\right)  \circ\Pi^{L}\left(  \overline{z},y\right)  \label{pp0}%
\end{equation}
and%
\begin{equation}
\Pi^{L}\left(  x,z\right)  \circ\Pi^{L}\left(  \overline{z},\overline
{x}\right)  =\Pi^{L}\left(  \overline{z},\overline{x}\right)  \circ\Pi
^{L}\left(  x,z\right)  =\operatorname*{id}\nolimits_{V}, \label{pp}%
\end{equation}
i.e. every $\Pi^{L}\left(  x,z\right)  $ is invertible and $\left(  \Pi
^{L}\left(  x,z\right)  \right)  ^{-1}=\Pi^{L}\left(  \overline{z}%
,\overline{x}\right)  $. This means that any left representation gives a
representation of a ternary group by a binary group.

Moreover, if a ternary group $(G,[\ ])$ is medial, then
\[
\Pi^{L}\left(  x_{1},x_{2}\right)  \circ\Pi^{L}\left(  x_{3},x_{4}\right)
=\Pi^{L}\left(  x_{3},x_{4}\right)  \circ\Pi^{L}\left(  x_{1},x_{2}\right)  ,
\]
i.e. obtained group is commutative. Indeed, by Corollary \ref{medial}, we
have
\begin{align*}
\Pi^{L}\left(  x_{1},x_{2}\right)  \circ\Pi^{L}\left(  x_{3},x_{4}\right)   &
=\Pi^{L}\left(  x_{1},x_{2}\right)  \circ\Pi^{L}\left(  x_{3},x_{4}\right)
\circ\Pi^{L}\left(  x,\overline{x}\right) \\
&  =\Pi^{L}\left(  \left[  [x_{1}x_{2}x_{3}]x_{4}\,x\right]  ,\,\overline
{x}\right)  =\Pi^{L}\left(  \left[  [x_{3}x_{4}x_{1}]x_{2}\,x\right]
,\,\overline{x}\right) \\
&  =\Pi^{L}\left(  x_{3},x_{4}\right)  \circ\Pi^{L}\left(  x_{1},x_{2}\right)
\circ\Pi^{L}\left(  x,\overline{x}\right) \\
&  =\Pi^{L}\left(  x_{3},x_{4}\right)  \circ\Pi^{L}\left(  x_{1},x_{2}\right)
.
\end{align*}

If $(G,[\ ])$ is commutative, then also $\Pi^{L}\left(  x,y\right)  =\Pi
^{L}\left(  y,x\right)  $, because
\begin{align*}
\Pi^{L}\left(  x,y\right)   &  =\Pi^{L}\left(  x,y\right)  \circ\Pi^{L}\left(
x,\overline{x}\right)  =\Pi^{L}\left(  \left[  x\,y\,x\right]  ,\,\overline
{x}\right) \\
&  =\Pi^{L}\left(  \left[  y\,x\,x\right]  ,\,\overline{x}\right)  =\Pi
^{L}\left(  y,x\right)  \circ\Pi^{L}\left(  x,\overline{x}\right)  =\Pi
^{L}\left(  y,x\right)  .
\end{align*}
Thus in the case of commutative and idempotent ternary groups any left
representation is idempotent and, in the consequence, $\left(  \Pi^{L}\left(
x,y\right)  \right)  ^{-1}=\Pi^{L}\left(  x,y\right)  $. This means that
commutative and idempotent ternary groups are represented by boolean groups.

\begin{proposition}
\label{xx}Let $\left(  G,\left[  \ \right]  \right)  =\mathrm{der\,}\left(
G,\odot\right)  $ be a ternary group derived from a binary group $\left(
G,\odot\right)  $. There is one-to-one correspondence between representations
of $\left(  G,\odot\right)  $ and left representations of $\left(  G,\left[
\ \right]  \right)  $.
\end{proposition}

\begin{proof}
Because $\left(  G,\left[  \;\right]  \right)  =\mathrm{der}\,\left(
G,\odot\right)  $, then $x\odot y=[xey]$ and $\overline{e}=e$, where $e$ is
unity of the binary group $\left(  G,\odot\right)  $. If $\pi\in
\mathrm{Rep}\left(  G,\odot\right)  $, then (as it is not difficult to see)
$\Pi^{L}\left(  x,y\right)  =\pi\left(  x\right)  \circ\pi\left(  y\right)  $
is a left representation of $\left(  G,\left[  \;\right]  \right)  $.
Conversely, if $\Pi^{L}$ is a left representation of $\left(  G,\left[
\;\right]  \right)  $ then $\pi\left(  x\right)  =\Pi^{L}\left(  x,e\right)  $
is a representation of $\left(  G,\odot\right)  $. Moreover, in this case
$\Pi^{L}\left(  x,y\right)  =\pi\left(  x\right)  \circ\pi\left(  y\right)  $.
Indeed, by Lemma \ref{xxl}, we have
\[
\Pi^{L}\left(  x,y\right)  =\Pi^{L}\left(  x,[eye]\right)  =\Pi^{L}\left(
[xey],e\right)  =\Pi^{L}\left(  x,e\right)  \circ\Pi^{L}\left(  y,e\right)
=\pi\left(  x\right)  \circ\pi\left(  y\right)
\]
for all $x,y\in G$.
\end{proof}

Let $(G,[\;])$ be a ternary group and $(G\times G,\ast)$ be a semigroup used
to the construction of left representations. According to Post \cite{pos} we
say that two pairs $(a,b)$, $(c,d)$ of elements of $G$ are equivalent, if
there exists an element $x\in G$ such that $[abx]=[cdx]$. Using a covering
group we can see that if this equation holds for some $x\in G$, then it holds
also for all $x\in G$. This means that
\[
\Pi^{L}(a,b)=\Pi^{L}(c,d)\Longleftrightarrow(a,b)\sim(c,d),
\]
%(vstav podhodiashchyi nomer formuly)
i.e.
\[
\Pi^{L}(a,b)=\Pi^{L}(c,d)\Longleftrightarrow\lbrack abx]=[cdx]
\]
for some $x\in G$. Indeed, if $[abx]=[cdx]$ holds for some $x\in G$, then
\begin{align*}
\Pi^{L}(a,b) &  =\Pi^{L}(a,b)\circ\Pi^{L}(x,\overline{x})=\Pi^{L}%
([abx],\,\overline{x})\\
&  =\Pi^{L}([cdx],\,\overline{x})=\Pi^{L}(c,d)\circ\Pi^{L}(x,\overline{x}%
)=\Pi^{L}(c,d).
\end{align*}
The converse is obvious.

Now we consider the second construction. Let $(G,[\ ])$ be a ternary group. On
$G\times G$ we define the following binary operation
\[
(x,y)\diamond(u,v)=(u,[vxy])
\]
Then $(G\times G,\diamond)$ is a binary semigroup which is isomorphic to
$(G\times G,\ast)$. This isomorphism has the form $\varphi((x,y))=\left(
\overline{y},\overline{x}\right)  $. Indeed,
\begin{align*}
\varphi((x,y)\diamond(u,v)) &  =\varphi(\,(u,[vxy])\,)=(\,\overline
{[vxy]},\overline{u}\,)\\
&  =(\,[\,\overline{y},\overline{x},\overline{v}\,],\overline{u}%
\,)=(\,\overline{y},\overline{x}\,)\ast(\,\overline{v},\overline{u}%
\,)=\varphi((x,y))\ast\varphi((u,v)).
\end{align*}

Basing on this construction we can define

\begin{definition}
A \textit{right representation} of a ternary group $G$ in $V$ is a map
$\Pi^{R}:G\times G\rightarrow\operatorname*{End}\,V$ such that%
\begin{align}
\Pi^{R}\left(  x_{3},x_{4}\right)  \circ\Pi^{R}\left(  x_{1},x_{2}\right)   &
=\Pi^{R}\left(  x_{1},\left[  x_{2}x_{3}x_{4}\right]  \right)  ,\;\;\forall
x_{1},x_{2},x_{3},x_{4}\in G\label{r1}\\
\Pi^{R}\left(  x,\overline{x}\right)   &  =\operatorname*{id}\nolimits_{V}%
,\;\;\forall x\in G. \label{r2}%
\end{align}

\end{definition}

From (\ref{r1})-(\ref{r2}) it follows that%
\begin{equation}
\Pi^{R}\left(  x,y\right)  =\Pi^{R}\left(  x,\left[  z\,\overline
{z}\,y\right]  \right)  =\Pi^{R}\left(  \overline{z},y\right)  \circ\Pi
^{R}\left(  x,z\right)  \ \ \forall\ x,y,z\in G. \label{pr}%
\end{equation}

It is easy to check that $\Pi^{R}\left(  x,y\right)  =\Pi^{L}\left(
\overline{y},\overline{x}\right)  =\left(  \Pi^{L}\left(  x,y\right)  \right)
^{-1}$. So it is enough to consider only left representations (as in binary case).

\begin{example}
\label{exam-kg}Let $G$ be a ternary group and $\mathbb{K}G$ is a vector space
spanned by $G$, which means that any element of $\mathbb{K}G$ can be uniquely
presented in the form $u=\sum_{i=1}^{n}k_{i}y_{i}$, $k_{i}\in\mathbb{K},$
$y_{i}\in G$. Then left and right regular representations are defined by
\begin{align}
\Pi_{reg}^{L}\left(  x_{1},x_{2}\right)  u  &  =\sum_{i=1}^{n}k_{i}\left[
x_{1}x_{2}y_{i}\right]  ,\label{pr1}\\
\Pi_{reg}^{R}\left(  x_{1},x_{2}\right)  u  &  =\sum_{i=1}^{n}k_{i}\left[
y_{i}x_{1}x_{2}\right]  , \label{pr2}%
\end{align}

\end{example}

\section{Middle representations}

Now we build another type of representations using the following construction.
For a given ternary group $(G,[\ ])$ we define on $G\times G^{op}$, where
$G^{op}$ is a ternary group having opposite multiplication, the following
ternary operation $\langle\ \rangle$ putting
\begin{equation}
\langle(x_{1},y_{1}),(x_{2},y_{2}),(x_{3},y_{3})\rangle=(\,[x_{1}%
\,x_{2}\,x_{3}],\,[y_{3}\,y_{2}\,y_{1}]\,)\label{mmm}%
\end{equation}
for all $x_{1},x_{2},x_{3},y_{1},y_{2},y_{3}\in G$. It is not difficult to see
that $(G\times G,\langle\ \rangle)$ is a ternary group as a direct product of
ternary groups. This group is commutative (medial, idempotent), if and only if
$(G,[\ ])$ is commutative (respectively: medial, idempotent). It is clear
that
\[
\langle(x,y),(\overline{x},\overline{y}),(a,b)\rangle=\langle
(a,b),(x,y),(\overline{x},\overline{y})\rangle=(a,b)
\]%
\[
\langle(\overline{x},\overline{y}),(x,y),(a,b)\rangle=\langle(a,b),(\overline
{x},\overline{y}),(x,y)\rangle=(a,b)
\]
for all $x,y,a,b\in G$. This means that in the group $(G\times G,\langle
\ \rangle)$ the element skew to $(x,y)$ has the form $(\overline{x}%
,\overline{y})$, where $\overline{x}$ is skew in $(G,[\ ])\,$.

Using (\ref{mmm}) we construct the middle representations as follows.

\begin{definition}
A \textit{middle representation} of a ternary group $G$ in $V$ is a map
$\Pi^{M}:G\times G\rightarrow\operatorname*{End}\,V$ such that%
\begin{align}
\Pi^{M}\left(  x_{3},y_{3}\right)  \circ\Pi^{M}\left(  x_{2},y_{2}\right)
\circ\Pi^{M}\left(  x_{1},y_{1}\right)   &  =\Pi^{M}\left(  \left[  x_{3}%
x_{2}x_{1}\right]  ,\left[  y_{1}y_{2}y_{3}\right]  \right)  ,\;\;\forall
x_{1},x_{2},x_{3},y_{1},y_{2},y_{3}\in G,\label{pm}\\
\Pi^{M}\left(  x,y\right)  \circ\Pi^{M}\left(  \overline{x},\overline
{y}\right)   &  =\Pi^{M}\left(  \overline{x},\overline{y}\right)  \circ\Pi
^{M}\left(  x,y\right)  =\operatorname*{id}\nolimits_{V}\ \ \forall\,x,y\in
G.\label{pm1}%
\end{align}

\end{definition}

It is seen that a middle representation is a ternary group homomorphism
$\Pi^{M}:G\times G^{op}\rightarrow\mathrm{der}\operatorname*{End}\,V$. Note
that instead of (\ref{pm1}) one can use $\Pi^{M}\left(  x,\overline{y}\right)
\circ\Pi^{M}\left(  \overline{x},y\right)  =id_{V}$ after changing $x$ to
$\overline{x}$ and taking into account that $x=\overline{\overline{x}}$.

\begin{remark}
In case elements $x$ and $y$ are idempotent we have $\Pi^{M}(x,y)\circ\Pi
^{M}(x,y)=id_{V}$, which means that the matrices $\Pi^{M}$ are Boolean. Thus
all middle representation matrices of idempotent ternary groups are Boolean.
\end{remark}

In general, the composition $\Pi^{M}\left(  x_{1},y_{1}\right)  \circ\Pi
^{M}\left(  x_{2},y_{2}\right)  $ is not a middle representation, but the
following proposition holds.

\begin{proposition}
If$\;\Pi^{M}$ is a middle representation of a ternary group $(G,[\ ])$, then
for any fixed $z\in G$

1. Let $\Pi_{z}^{L}(x,y)=\Pi^{M}(x,z)\circ\Pi^{M}(y,\overline{z})$ is a left
representation of $(G,[\ ])$, then $\Pi_{z}^{L}(x,y)\circ\Pi_{z^{\prime}}%
^{L}(x^{\prime},y^{\prime})=\Pi_{z^{\prime}}^{L}(\left[  xyz^{\prime}\right]
,y^{\prime})$.

2. Let $\Pi_{z}^{R}(x,y)=\Pi^{M}(z,y)\circ\Pi^{M}(\overline{z},x)$ is a right
representation of $(G,[\ ])$, then $\Pi_{z}^{R}(x,y)\circ\Pi_{z^{\prime}}%
^{R}(x^{\prime},y^{\prime})=\Pi_{z}^{R}(x,\left[  yx^{\prime}y^{\prime
}\right]  )$.
\end{proposition}

\begin{proof}
The proof is a verification of the corresponding axioms.
\end{proof}

In particular, $\Pi_{z}^{L}$ ($\Pi_{z}^{R}$) is a family of left (right) representations.

\begin{corollary}
If a middle representation $\Pi^{M}$ of a ternary group $(G,[\ ])$ satisfies
$\Pi^{M}\left(  x,\overline{x}\right)  =\mathrm{id}_{V}$ for all $x\in G$,
then it is a left and right representation and $\Pi^{M}(x,y)=\Pi^{M}(y,x)$ for
all $x,y\in G$.
\end{corollary}

\begin{proof}
Indeed,
\begin{align*}
\Pi^{M}(x,y) &  =\Pi^{M}([xy\overline{y}],[y\overline{z}z])\\
&  =\Pi^{M}(x,z)\circ\Pi^{M}(y,\overline{z})\circ\Pi^{M}(\overline{y}%
,y)=\Pi^{M}(x,z)\circ\Pi^{M}(y,\overline{z})=L(x,y).
\end{align*}
Similarly
\begin{align*}
\Pi^{M}(x,y) &  =\Pi^{M}([z\overline{z}x],[\overline{x}xy])\\
&  =\Pi^{M}(z,y)\circ\Pi^{M}(\overline{z},x)\circ\Pi^{M}(x,\overline{x}%
)=\Pi^{M}(z,y)\circ\Pi^{M}(\overline{z},x)=R(x,y)
\end{align*}
and
\begin{align*}
\Pi^{M}(x,y) &  =\Pi^{M}(\,[\,x\,y\,\overline{y}\,],[\,y\,x\,\overline
{x}\,]\,)\\
&  =\Pi^{M}(x,\overline{x})\circ\Pi^{M}(y,x)\circ\Pi^{M}(\overline{y}%
,y)=\Pi^{M}(y,x),
\end{align*}
which completes the proof.
\end{proof}

Observe that in general $\Pi_{reg}^{M}(x,\overline{x})\neq\mathrm{id}$. 

It can be shown that for regular representations we have the following
commutation relations%
\[
\Pi_{reg}^{L}\left(  x_{1},y_{1}\right)  \circ\Pi_{reg}^{R}\left(  x_{2}%
,y_{2}\right)  =\Pi_{reg}^{R}\left(  x_{2},y_{2}\right)  \circ\Pi_{reg}%
^{L}\left(  x_{1},y_{1}\right)  .
\]

\begin{proposition}
For a finite (or countable) ternary group $\left(  G,\left[  \;\right]
\right)  $ left and right representations are unitary.
\end{proposition}

\begin{proof}
Take a scalar product $\left\langle \;,\;\right\rangle $ in $\mathbb{K}G$
which makes $G$ an orthonormal basis, i.e. $\left\langle g,h\right\rangle
=\delta_{g,h}$. Then the unitarity follows from uniqueness of solutions to the
group equations $\left[  xyg\right]  =h$ (see (\ref{ac})).
\end{proof}

\section{Relation between representations}

Let $(G,[\;])$ be a ternary group and let $(G\times G,\langle\;\rangle)$ be a
ternary group used to the construction of the middle representation. In
$(G,[\;])$ (and in the consequence in $(G\times G,\langle\;\rangle)\;$) we
define the relation
\[
(a,b)\sim(c,d)\Longleftrightarrow\lbrack azb]=[czd]
\]
for all $z\in G$. It is not difficult to see that this relation is a
congruence in $(G\times G,\langle\;\rangle)$. For regular representations
$\Pi_{reg}^{M}(a,b)=\Pi_{reg}^{M}(c,d)$ if $(a,b)\sim(c,d)$.

Thus in Example \ref{exam-reg} we have $\Pi^{M}(a,b)=\Pi^{M}%
(c,d)\Longleftrightarrow(a+b)=(c+d)\operatorname{mod}\,3$. Hence, the
computation of middle representations can be reduced to the computation only
of three cases $\Pi^{M}(0,0)$, $\Pi^{M}(1,0)$, $\Pi^{M}(2,0)$.

So we have the following relation
\[
a\eqsim a^{\prime}\Longleftrightarrow a=[\overline{x}a^{\prime}%
x]\;\;\text{for\ some}\;x\in G
\]
or equivalently
\[
a\eqsim a^{\prime}\Longleftrightarrow a^{\prime}=[xa\overline{x}%
]\;\;\text{for\ some}\;x\in G.
\]

It is not difficult to see that it is an equivalence relation on $(G,[\;])$,
moreover, if $(G,[\;])$ is medial, then this relation is a congruence.

Let $(G\times G,\langle\;\rangle)$ be a ternary group used for a construction
of middle representations, then
\begin{align*}
(a,b) &  \eqsim(a^{\prime},b)\Longleftrightarrow a^{\prime}=[xa\overline
{x}]\;\;\text{and}\;\;\\
b^{\prime} &  =[yb\overline{b}\,]\;\;\text{for\ some}\;(x,y)\in G\times G
\end{align*}
is an equivalence relation on $(G\times G,\langle\;\rangle)$. Moreover, if
$(G,[\;])$ is medial, then this relation is a congruence. Unfortunately, it is
a weak relation. In a ternary group $\mathbb{Z}_{3}$, where $[xyz]=\left(
x-y+z\right)  (\operatorname{mod}3)$ we have only one class, i.e. all elements
are equivalent. In $\mathbb{Z}_{4}$ with the operation $[xyz]=\left(
x+y+z+1\right)  (\operatorname{mod}\,4)$ we have $a\eqsim a^{\prime
}\Longleftrightarrow a=a^{\prime}$. But for this relation holds the following

\begin{lemma}
If $(a,b)\eqsim(a^{\prime},b^{\prime})$, then
\[
\operatorname{tr}\Pi^{M}(a,b)=\operatorname{tr}\Pi^{M}(a^{\prime},b^{\prime}).
\]

\end{lemma}

\begin{proof}
We have $\operatorname{tr}(AB)=\operatorname{tr}(BA)$ for all $A,B\in
\mathrm{End}V$. Indeed,
\begin{align*}
\operatorname{tr}\Pi^{M}(a,b)  &  =\operatorname{tr}\Pi^{M}([xa^{\prime
}\overline{x}],[yb^{\prime}\overline{y}]\,)=\operatorname{tr}\left(  \Pi
^{M}(x,\overline{y})\circ\Pi^{M}(a^{\prime},b^{\prime})\circ\Pi^{M}%
(\overline{x},y)\right) \\
&  =\operatorname{tr}\left(  \Pi^{M}(x,\overline{y})\circ\Pi^{M}(\overline
{x},y)\circ\Pi^{M}(a^{\prime},b^{\prime})\right)  =\operatorname{tr}\left(
id_{V}\circ\Pi^{M}(a^{\prime}b^{\prime})\right) \\
&  =\operatorname{tr}\Pi^{M}(a^{\prime},b^{\prime})
\end{align*}

\end{proof}

We can \textquotedblleft algebralize\textquotedblright\ the above regular
representations from the Example \ref{exam-reg} in the following way. From
(\ref{p1}) we have for the left representation $\Pi_{reg}^{L}\left(
i,j\right)  \circ\Pi_{reg}^{L}\left(  k,l\right)  =\Pi_{reg}^{L}\left(
i,\left[  jkl\right]  \right)  $, where $\left[  jkl\right]  =j-k+l$,
$i,j,k,l\in\mathbb{Z}_{3}$. Denote $\gamma_{i}^{L}=\Pi_{reg}^{L}\left(
0,i\right)  $, $i\in\mathbb{Z}_{3}$, then we obtain the algebra with the
relations $\gamma_{i}^{L}\gamma_{j}^{L}=\gamma_{i+j}^{L}$. Conversely, any
matrix representation of $\gamma_{i}\gamma_{j}=\gamma_{i+j}$ leads to the left
representation by $\Pi^{L}\left(  i,j\right)  =\gamma_{j-i}$.

In the case of the middle regular representation we introduce $\gamma
_{k+l}^{M}=\Pi_{reg}^{M}\left(  k,l\right)  $, $k,l\in\mathbb{Z}_{3}$, then we
obtain
\begin{equation}
\gamma_{i}^{M}\gamma_{j}^{M}\gamma_{k}^{M}=\gamma_{\left[  ijk\right]  }%
^{M},\;\;\;i,j,k\in\mathbb{Z}_{3}. \label{ggg}%
\end{equation}

In some sense (\ref{ggg}) can be treated as a \textit{ternary analog of
Clifford algebra}. As before, any matrix representation of (\ref{ggg}) gives
the middle representation $\Pi^{M}\left(  k,l\right)  =\gamma_{k+l}$.

In our \textit{derived} case the connection with the standard group
representations is given by

\begin{proposition}
Let $\left(  G,\odot\right)  $ be a binary group, and the ternary derived
group as $\left(  G,\left[  \;\right]  \right)  =\mathrm{der}\,\left(
G,\odot\right)  $. There is one-to-one correspondence between a pair of
commuting binary groups representations and a middle ternary derived group representation.
\end{proposition}

\begin{proof}
Let $\pi,\rho\in Rep\left(  G,\odot\right)  $, $\pi\left(  x\right)  \circ
\rho\left(  y\right)  =\rho\left(  y\right)  \circ\pi\left(  x\right)  $ and
$\Pi^{L}\in Rep\left(  G,\left[  \;\right]  \right)  $. We take
\begin{align*}
\Pi^{M}\left(  x,y\right)   &  =\pi\left(  x\right)  \circ\rho\left(
y^{-1}\right)  ,\\
\pi\left(  x\right)   &  =\Pi^{M}\left(  x,e\right)  ,\\
\rho\left(  x\right)   &  =\Pi^{M}\left(  e,\overline{x}\right)  .
\end{align*}
Then using (\ref{pm}) we prove the needed representation laws.
\end{proof}

Let $(G,[\;])$ be a fixed ternary group, $(G\times G,\langle\;\rangle)$ a
corresponding ternary group used in the construction of middle
representations, $(\left(  G\times G\right)  ^{\ast},\circledast)$ a covering
group of $(G\times G,\langle\;\rangle)$, $(G\times G,\diamond)=$\textrm{$ret$%
}$_{(a,b)}(G\times G,\langle\;\rangle)$. If $\Pi^{M}(a,b)$ is a middle
representation of $(G,[\;])$, then $\pi$ defined by
\[
\pi(x,y,0)=\Pi^{M}(x,y)
\]
and
\[
\pi(x,y,1)=\Pi^{M}(x,y)\circ\Pi^{M}(a,b)
\]
is a representation of the covering group. Moreover
\[
\rho(x,y)=\Pi^{M}(x,y)\circ\Pi^{M}(a,b)=\pi(x,y,1)
\]
is a representation of the above retract induced by $(a,b)$. Indeed,
$(\overline{a},\overline{b})$ is the identity of this retract and
$\rho(\overline{a},\overline{b})=\Pi^{M}(\overline{a},\overline{b})\circ
\Pi^{M}(a,b)=\operatorname*{id}\nolimits_{V}$. Similarly
\begin{align*}
\rho\left(  (x,y)\diamond(z,u)\right)   &  =\rho\left(  \langle
(x,y),(a,b),(z,u)\rangle\right)  =\rho\left(  \lbrack xaz],[uby]\right)  \\
&  =\Pi^{M}\left(  [xaz],[uby]\right)  )\circ\Pi^{M}(a,b)\\
&  =\Pi^{M}(x,y)\circ\Pi^{M}(a,b)\circ\Pi^{M}(z,u)\circ\Pi^{M}(a,b)\\
&  =\rho(x,y)\circ\rho(z,u)
\end{align*}

But $\tau(x)=(x,\overline{x})$ is an embedding of $(G,[\;])$ into $(G\times
G,\langle\;\rangle)$. Hence $\mu$ defined by $\mu(x,0)=\Pi^{M}(x,\overline
{x})$ and $\mu(x,1)=\Pi^{M}(x,\overline{x})\circ\Pi^{M}(a,\overline{a})$ is a
representation of a covering group $G^{\ast}$ for $(G,[\;])$ (see Post theorem
for $a=c$). On the other hand, $\beta(x)=\Pi^{M}(x,\overline{x})\circ\Pi
^{M}(a,\overline{a})$ is a representation of a binary retract $(G,\cdot
\,)=$\textrm{$ret$}$_{a}(G,[\;])$. That $\beta$ can induce some middle
representation of $(G,[\;])$ (by Gluskin-Hossz\'{u} theorem).

Note that in a ternary group of quaternions $(\mathbb{K},[\;])$, where
$[xyz]=xyz(-1)=-xyz$ and $xy$ is a multiplication of quaternions ($-1$ is a
central element) we have $\overline{1}=-1$, $\overline{-1}=1$ and
$\overline{x}=x$ for others. In $(K\times K,\langle\;\rangle)$ we have
$(a,b)\sim(-a,-b)$ and $(a,-b)\sim(-a,b)$, which gives 32 two-elements
equivalence classes. The embedding $\tau(x)=(x,\overline{x})$ suggest that
$\Pi^{M}(i,i)=\pi(i)\neq\pi(-i)=\Pi^{M}(-i,-i)$. Generally $\Pi^{M}%
(a,b)\neq\Pi^{M}(-a,-b)$ and $\Pi^{M}(a,-b)\neq\Pi^{M}(-a,b)$.

The relation $(a,b)\sim(c,d)\Longleftrightarrow\lbrack abx]=[cdx]$ for all
$x\in G$ is a congruence on $(G\times G,\ast)$. Note that this relation can be
defined as "for some $x$". Indeed, using a covering group we can see that if
$[abx]=[cdx]$ holds for some $x$ then holds also for all $x$. Thus $\pi
^{L}(a,b)=\Pi^{L}(c,d)\Longleftrightarrow(a,b)\sim(c,d)$. Indeed
\begin{align*}
\Pi^{L}(a,b)  &  =\Pi^{L}(a,b)\circ\Pi^{L}(x,\overline{x})=\Pi^{L}%
([a\ b\ x],\overline{x}))\\
&  =\Pi^{L}([c\ d\ x],\overline{x}))=\Pi^{L}(c,d)\circ\Pi^{L}(x,\overline
{x})=\Pi^{L}(c,d).
\end{align*}

\begin{proposition}
Every left representation of a commutative group $(G,[\;])$ is a middle representation.
\end{proposition}

\begin{proof}
Indeed,
\begin{align*}
\Pi^{L}(x,y)\circ\Pi^{L}(\overline{x},\overline{y})  &  =\Pi^{L}%
([x\ y\ \overline{x}],\overline{y})\\
&  =\Pi^{L}([x\ \overline{x}\ y],\overline{y})=\Pi^{L}(y,\overline
{y})=\mathrm{id}_{V}%
\end{align*}
and
\begin{align*}
\Pi^{L}(x_{1},x_{2})\circ\Pi^{L}(x_{3},x_{4})\circ\Pi^{L}(x_{5},x_{6})  &
=\Pi^{L}([[x_{1}x_{2}x_{3}]x_{4}x_{5}],x_{6})\\
&  =\Pi^{L}([[x_{1}x_{3}x_{2}]x_{4}x_{5}],x_{6})=\Pi^{L}([x_{1}x_{3}%
[x_{2}x_{4}x_{5}]],x_{6})
\end{align*}%
\[
=\Pi^{L}([x_{1}x_{3}[x_{5}x_{4}x_{2}]],x_{6})=\Pi^{L}([x_{1}x_{3}x_{5}%
],[x_{4}x_{2}x_{6}])=\Pi^{L}([x_{1}x_{3}x_{5}],[x_{6}x_{4}x_{2}]).
\]

\end{proof}

Note that the converse holds only for middle representations such that
$\Pi^{M}(x,\overline{x})=\mathrm{id}_{V}$.

\begin{theorem}
There is one-one-correspondence between left representations of $(G,[\;])$ and
binary representations of the retract \textrm{$ret$}$_{a}(G,[\;])$.
\end{theorem}

\begin{proof}
Let $\Pi^{L}(x,a)$ is given, then define $\rho(x)=\Pi^{L}(x,a)$ is such
representation of the retract which can be directly shown. Conversely, assume
that $\rho(x)$ is a representation of the retract \textrm{$ret$}$_{a}%
(G,[\;])$. Define $\Pi^{L}(x,y)=\rho(x)\circ\rho(\overline{y})^{-1}$, then
$\Pi^{L}(x,y)\circ\Pi^{L}(z,u)=\rho(x)\circ\rho(\overline{y})^{-1}\circ
\rho(z)\circ\rho(\overline{u})^{-1}=\rho(x\circledast(\overline{y})^{-1}%
\circ\circledast z)\circ\rho(\overline{u})^{-1}=\rho([\,[\,x\,a\,[\,\overline
{a}\,y\,\overline{a}\,]\,]\,a\,z\,])\circ\rho(\overline{u})^{-1}%
=\rho([\,x\,y\,x\,])\circ\rho(\overline{u})^{-1}=\Pi^{L}([\,x\,y\,z\,],u)$
which completes the proof.
\end{proof}

\begin{remark}
It is seen that Proposition \ref{xx} is a direct consequence of this theorem.
\end{remark}

\section{Matrix representations}

Now we give examples of matrix representations for concrete ternary groups.

\begin{example}
\label{exam-reg}Let $G=\mathbb{Z}_{3}\ni\left\{  0,1,2\right\}  $ and the
ternary multiplication is $\left[  xyz\right]  =x-y+z$. Then $\left[
xyz\right]  =\left[  zyx\right]  $ and $\overline{0}=0,$ $\overline{1}=1,$
$\overline{2}=2$, therefore $(G,[\ ])$ is an idempotent medial ternary group.
Thus $\Pi^{L}(x,y)=\Pi^{R}(y,x)$ and
\begin{equation}
\Pi^{L}(a,b)=\Pi^{L}(c,d)\Longleftrightarrow
(a-b)=(c-d)\mathrm{\operatorname{mod}}\,3.\label{mod3}%
\end{equation}
Straightforward calculations give the left regular representation in the
manifest matrix form%
\begin{align*}
\Pi_{reg}^{L}\left(  0,0\right)   &  =\Pi_{reg}^{L}\left(  2,2\right)
=\Pi_{reg}^{L}\left(  1,1\right)  =\Pi_{reg}^{R}\left(  0,0\right)  \\
&  =\Pi_{reg}^{R}\left(  2,2\right)  =\Pi_{reg}^{R}\left(  1,1\right)
=\left(
\begin{array}
[c]{ccc}%
1 & 0 & 0\\
0 & 1 & 0\\
0 & 0 & 1
\end{array}
\right)  \\
&  =[1]\oplus\lbrack1]\oplus\lbrack1],\\
\Pi_{reg}^{L}\left(  2,0\right)   &  =\Pi_{reg}^{L}\left(  1,2\right)
=\Pi_{reg}^{L}\left(  0,1\right)  =\Pi_{reg}^{R}\left(  2,1\right)  \\
&  =\Pi_{reg}^{R}\left(  1,0\right)  =\Pi_{reg}^{R}\left(  0,2\right)
=\left(
\begin{array}
[c]{ccc}%
0 & 1 & 0\\
0 & 0 & 1\\
1 & 0 & 0
\end{array}
\right)  \\
&  =[1]\oplus\left(
\begin{array}
[c]{cc}%
-\dfrac{1}{2} & -\dfrac{\sqrt{3}}{2}\\
\dfrac{\sqrt{3}}{2} & -\dfrac{1}{2}%
\end{array}
\right)  =[1]\oplus\left[  -\dfrac{1}{2}+\dfrac{1}{2}i\sqrt{3}\right]
\oplus\left[  -\dfrac{1}{2}-\dfrac{1}{2}i\sqrt{3}\right]  ,\\
\Pi_{reg}^{L}\left(  2,1\right)   &  =\Pi_{reg}^{L}\left(  1,0\right)
=\Pi_{reg}^{L}\left(  0,2\right)  =\Pi_{reg}^{R}\left(  2,0\right)  \\
&  =\Pi_{reg}^{R}\left(  1,2\right)  =\Pi_{reg}^{R}\left(  0,1\right)
=\left(
\begin{array}
[c]{ccc}%
0 & 0 & 1\\
1 & 0 & 0\\
0 & 1 & 0
\end{array}
\right)  \\
&  =[1]\oplus\left(
\begin{array}
[c]{cc}%
-\dfrac{1}{2} & \dfrac{\sqrt{3}}{2}\\
-\dfrac{\sqrt{3}}{2} & -\dfrac{1}{2}%
\end{array}
\right)  =[1]\oplus\left[  -\dfrac{1}{2}-\dfrac{1}{2}i\sqrt{3}\right]
\oplus\left[  -\dfrac{1}{2}+\dfrac{1}{2}i\sqrt{3}\right]  .
\end{align*}

\end{example}

Consider the middle representation constructions of Examples \ref{exam-kg} and
\ref{exam-reg}.

\begin{example}
The middle regular representations from Example \ref{exam-kg} is defined by
\[
\Pi_{reg}^{M}\left(  x_{1},x_{2}\right)  u=\sum_{i=1}^{n}k_{i}\left[
x_{1}y_{i}x_{2}\right]
\]
For regular representations we have
\begin{align}
\Pi_{reg}^{M}\left(  x_{1},y_{1}\right)  \circ\Pi_{reg}^{R}\left(  x_{2}%
,y_{2}\right)   &  =\Pi_{reg}^{R}\left(  y_{2},y_{1}\right)  \circ\Pi
_{reg}^{M}\left(  x_{1},x_{2}\right)  ,\label{pp30}\\
\Pi_{reg}^{M}\left(  x_{1},y_{1}\right)  \circ\Pi_{reg}^{L}\left(  x_{2}%
,y_{2}\right)   &  =\Pi_{reg}^{L}\left(  x_{1},x_{2}\right)  \circ\Pi
_{reg}^{M}\left(  y_{2},y_{1}\right)  . \label{pp3}%
\end{align}

\end{example}

\begin{example}
For the middle regular representation matrices we obtain%
\begin{align*}
\Pi_{reg}^{M}\left(  0,0\right)   &  =\Pi_{reg}^{M}\left(  1,2\right)
=\Pi_{reg}^{M}\left(  2,1\right)  =\left(
\begin{array}
[c]{ccc}%
1 & 0 & 0\\
0 & 0 & 1\\
0 & 1 & 0
\end{array}
\right)  ,\\
\Pi_{reg}^{M}\left(  0,1\right)   &  =\Pi_{reg}^{M}\left(  1,0\right)
=\Pi_{reg}^{M}\left(  2,2\right)  =\left(
\begin{array}
[c]{ccc}%
0 & 1 & 0\\
1 & 0 & 0\\
0 & 0 & 1
\end{array}
\right)  ,\\
\Pi_{reg}^{M}\left(  0,2\right)   &  =\Pi_{reg}^{M}\left(  2,0\right)
=\Pi_{reg}^{M}\left(  1,1\right)  =\left(
\begin{array}
[c]{ccc}%
0 & 0 & 1\\
0 & 1 & 0\\
1 & 0 & 0
\end{array}
\right)  .
\end{align*}

\end{example}

\begin{example}
The above representation $\Pi_{reg}^{M}$ of $\left(  \mathbb{Z}_{3},\left[
\;\right]  \right)  $ is equivalent to the orthogonal direct sum of two
irreducible representations%
\begin{align*}
\Pi_{reg}^{M}\left(  0,0\right)   &  =\Pi_{reg}^{M}\left(  1,2\right)
=\Pi_{reg}^{M}\left(  2,1\right)  =\left[  1\right]  \oplus\left[
\begin{array}
[c]{cc}%
-1 & 0\\
0 & 1
\end{array}
\right]  ,\\
\Pi_{reg}^{M}\left(  0,1\right)   &  =\Pi_{reg}^{M}\left(  1,0\right)
=\Pi_{reg}^{M}\left(  2,2\right)  =\left[  1\right]  \oplus\left[
\begin{array}
[c]{cc}%
\dfrac{1}{2} & -\dfrac{\sqrt{3}}{2}\\
-\dfrac{\sqrt{3}}{2} & -\dfrac{1}{2}%
\end{array}
\right]  ,\\
\Pi_{reg}^{M}\left(  0,2\right)   &  =\Pi_{reg}^{M}\left(  2,0\right)
=\Pi_{reg}^{M}\left(  1,1\right)  =\left[  1\right]  \oplus\left[
\begin{array}
[c]{cc}%
\dfrac{1}{2} & \dfrac{\sqrt{3}}{2}\\
\dfrac{\sqrt{3}}{2} & -\dfrac{1}{2}%
\end{array}
\right]  ,
\end{align*}
i.e. one-dimensional trivial $\left[  1\right]  $ and two-dimensional irreducible.
\end{example}

\begin{remark}
In this example $\Pi^{M}(x,\overline{x})=\Pi^{M}(x,x)\neq id_{V}$, but
$\Pi^{M}(x,y)\circ\Pi^{M}(x,y)=id_{V}$, and so $\Pi^{M}$ are of second degree.
\end{remark}

Let us consider a more complicated example of left representations.

\begin{example}
\label{exam-z4}Let $G=\mathbb{Z}_{4}\ni\left\{  0,1,2,3\right\}  $ and the
ternary multiplication is
\begin{equation}
\left[  xyz\right]  =\left(  x+y+z+1\right)  \operatorname{mod}4.
\label{mult4}%
\end{equation}
We have the multiplication table
\begin{align*}
\left[  x,y,0\right]   &  =\left(
\begin{array}
[c]{cccc}%
1 & 2 & 3 & 0\\
2 & 3 & 0 & 1\\
3 & 0 & 1 & 2\\
0 & 1 & 2 & 3
\end{array}
\right)  \ \ \ \ \ \ \ \left[  x,y,1\right]  =\left(
\begin{array}
[c]{cccc}%
2 & 3 & 0 & 1\\
3 & 0 & 1 & 2\\
0 & 1 & 2 & 3\\
1 & 2 & 3 & 0
\end{array}
\right) \\
\left[  x,y,2\right]   &  =\left(
\begin{array}
[c]{cccc}%
3 & 0 & 1 & 2\\
0 & 1 & 2 & 3\\
1 & 2 & 3 & 0\\
2 & 3 & 0 & 1
\end{array}
\right)  \ \ \ \ \ \ \ \left[  x,y,3\right]  =\left(
\begin{array}
[c]{cccc}%
0 & 1 & 2 & 3\\
1 & 2 & 3 & 0\\
2 & 3 & 0 & 1\\
3 & 0 & 1 & 2
\end{array}
\right)
\end{align*}
Then the skew elements are $\overline{0}=3,$ $\overline{1}=2,$ $\overline
{2}=1,$ $\overline{3}=0$, therefore $(G,[\ ])$ is an (nonidempotent)
commutative ternary group. The left representation is defined by expansion
$\Pi_{reg}^{L}\left(  x_{1},x_{2}\right)  u=\sum_{i=1}^{n}k_{i}\left[
x_{1}x_{2}y_{i}\right]  $, which means that
\[
\Pi_{reg}^{L}\left(  x,y\right)  |z>=|\left[  xyz\right]  >.
\]
Analogously, for right and middle representations
\[
\Pi_{reg}^{R}\left(  x,y\right)  |z>=|\left[  zxy\right]  >,\ \ \ \Pi
_{reg}^{M}\left(  x,y\right)  |z>=|\left[  xzy\right]  >.
\]
Therefore $|\left[  xyz\right]  >=|\left[  zxy\right]  >=|\left[  xzy\right]
>$ and
\[
\Pi_{reg}^{L}\left(  x,y\right)  =\Pi_{reg}^{R}\left(  x,y\right)
|z>=\Pi_{reg}^{M}\left(  x,y\right)  |z>,
\]
so $\Pi_{reg}^{L}\left(  x,y\right)  =\Pi_{reg}^{R}\left(  x,y\right)
=\Pi_{reg}^{M}\left(  x,y\right)  $. Thus it is sufficient to consider the
left representation only.

In this case the equivalence is $\Pi^{L}(a,b)=\Pi^{L}(c,d)\Longleftrightarrow
(a+b)=(c+d)$\textrm{$\operatorname{mod}$}$\,4$, and we obtain the following
classes%
\begin{align*}
\Pi_{reg}^{L}\left(  0,0\right)   &  =\Pi_{reg}^{L}\left(  1,3\right)
=\Pi_{reg}^{L}\left(  2,2\right)  =\Pi_{reg}^{L}\left(  3,1\right)  \\
&  =\left(
\begin{array}
[c]{cccc}%
0 & 0 & 0 & 1\\
1 & 0 & 0 & 0\\
0 & 1 & 0 & 0\\
0 & 0 & 1 & 0
\end{array}
\right)  =\left[  1\right]  \oplus\left[  -1\right]  \oplus\left[  -i\right]
\oplus\left[  i\right]  ,\\
\Pi_{reg}^{L}\left(  0,1\right)   &  =\Pi_{reg}^{L}\left(  1,0\right)
=\Pi_{reg}^{L}\left(  2,3\right)  =\Pi_{reg}^{L}\left(  3,2\right)  \\
&  =\left(
\begin{array}
[c]{cccc}%
0 & 0 & 1 & 0\\
0 & 0 & 0 & 1\\
1 & 0 & 0 & 0\\
0 & 1 & 0 & 0
\end{array}
\right)  =\left[  1\right]  \oplus\left[  -1\right]  \oplus\left[  -1\right]
\oplus\left[  -1\right]  ,\\
\Pi_{reg}^{L}\left(  0,2\right)   &  =\Pi_{reg}^{L}\left(  1,1\right)
=\Pi_{reg}^{L}\left(  2,0\right)  =\Pi_{reg}^{L}\left(  3,3\right)  \\
&  =\left(
\begin{array}
[c]{cccc}%
0 & 1 & 0 & 0\\
0 & 0 & 1 & 0\\
0 & 0 & 0 & 1\\
1 & 0 & 0 & 0
\end{array}
\right)  =\left[  1\right]  \oplus\left[  -1\right]  \oplus\left[  i\right]
\oplus\left[  -i\right]  ,\\
\Pi_{reg}^{L}\left(  0,3\right)   &  =\Pi_{reg}^{L}\left(  1,2\right)
=\Pi_{reg}^{L}\left(  2,1\right)  =\Pi_{reg}^{L}\left(  3,0\right)  \\
&  =\left(
\begin{array}
[c]{cccc}%
1 & 0 & 0 & 0\\
0 & 1 & 0 & 0\\
0 & 0 & 1 & 0\\
0 & 0 & 0 & 1
\end{array}
\right)  =\left[  1\right]  \oplus\left[  -1\right]  \oplus\left[  1\right]
\oplus\left[  1\right]  .
\end{align*}
It is seen that, due to the fact that the ternary operation (\ref{mult4}) is
commutative, there are only one-dimensional irreducible left representations.
\end{example}

\medskip In a similar way one can extend other notions of the classical group
representation theory to the ternary group case. This includes, e.g. direct
sum and tensor product of representations, characters, irreducibility (Schur
lemma), equivalence of representations etc.

\textbf{Acknowledgments}. A.B. is grateful to K. G\l azek and Z. Oziewicz for
interesting discussions, and S.D. would like to thank Jerzy Lukierski for kind
hospitality at the University of Wroc\l aw.

\end{document}